\documentclass[12pt]{article}
\usepackage{amssymb}
\usepackage{amsfonts}
\usepackage{graphicx}
\usepackage{amsmath}
\usepackage{booktabs}
\usepackage{multirow}
\newtheorem{theorem}{Theorem}

\newtheorem{algorithm}[theorem]{Algorithm}

\newtheorem{definition}[theorem]{Definition}

\newtheorem{proposition}[theorem]{Proposition}

\newenvironment{proof}[1][Proof]{\noindent\textbf{#1.} }{\ \rule{0.5em}{0.5em}}

\begin{document}
\title{Some notes on the trapezoidal rule for Fourier type integrals}
\author{Eleonora Denich \thanks{Dipartimento di Matematica e Geoscienze, Universit\`{a} di Trieste, Trieste, Italy, eleonora.denich@phd.units.it} \and Paolo Novati \thanks{Dipartimento di Matematica e Geoscienze, Universit\`{a} di Trieste, Trieste, Italy, novati@units.it}}
\date{}
\maketitle

\begin{abstract}
    This paper deals with the error analysis of the trapezoidal rule for the computation of Fourier type integrals, based on two double exponential transformations.
    The theory allows to construct algorithms in which the steplength and the number of nodes can be a priori selected.
    The analysis is also used to design an automatic integrator that can be employed without any knowledge of the function involved in the problem.
    Several numerical examples, which confirm the reliability of this strategy, are reported.
    
\end{abstract}

\section{Introduction}

In this work we consider the computation of the cosine and sine transform, defined respectively by
\begin{align*}
    F^{(c)}(\omega) &= \int_0^{+\infty} f(x) \cos \left( \omega x \right), \\
    F^{(s)}(\omega) &= \int_0^{+\infty} f(x) \sin \left( \omega x \right),
\end{align*}
where $f$ is an integrable function and $\omega \in \mathbb{R}$ (for an overview see e.g. \cite{Y} and the reference therein).
It is well known that evaluating these kind of integrals by means of standard quadrature techniques may lead to quite inaccurate results, especially when the function $f$ exhibits a slow decay and/or when the frequency $\omega$ is rather large.
For this reason, Ooura and Mori in \cite{OM,OM1,OM2} introduced two special exponential type transformations $x = \frac{\tau}{\omega} \phi_i(\xi)$, $i=1,2$, $\phi_i \colon \left(-\infty, +\infty \right) \rightarrow \left( 0, +\infty \right)$, $\tau >0$, defined by
\begin{align*}
    \phi_1(\xi) &= \frac{\xi}{1-e^{- 2 \pi \sinh \xi}}, \\
    \phi_2(\xi) &= \frac{\xi}{1-e^{-2 \xi-\alpha \left( 1-e^{-t} \right) - \beta \left( e^t-1 \right)}}, \quad 0<\alpha<\beta<1. 
\end{align*}
By using these transformations, the above integrals can be efficiently computed by employing the truncated trapezoidal rule in the following way:
\begin{align} 
    F^{(c)}(\omega) &\approx \frac{\tau}{\omega} h \sum_{j = -M}^{N} f \left( \frac{\tau}{\omega} \phi_i\left(jh -\frac{\pi}{2 \tau} \right) \right) \cos \left( \tau \phi_i\left(jh -\frac{\pi}{2 \tau} \right) \right) \phi_i^{'}\left(jh -\frac{\pi}{2 \tau} \right),\label{T} \\
    F^{(s)}(\omega) &\approx \frac{\tau}{\omega} h \sum_{j=-M}^N f \left( \frac{\tau}{\omega} \phi_i (jh) \right) \sin \left(\tau \phi_i(jh) \right) \phi_i^{'}(jh). \label{1 bis}
\end{align}
By considering a generalization of $\phi_2$, in \cite{O} Ooura proposed a method for the computation of the Fourier transform
\begin{equation*}
    F(\omega) = \int_0^{+ \infty} f(x) e^{i \omega x} dx = F^{(c)}(\omega)+ i F^{(s)}(\omega),
\end{equation*}
where the function $f$ may have singularities or exhibits slow convergence at infinity.
For an overview of the most commonly employed techniques for integrals involving oscillating functions we quote here \cite{MS,AH,Wong} and the references therein.

In this framework, here we present reliable error bounds that turn out from a proper selection of the parameters $h,M,N$, that allows to equalize the error contributions arising from formulas (\ref{T})-(\ref{1 bis}).
The approximations presented require the knowledge of the region of analiticity of the function $f$.
After showing the standard approach for the error analysis as given in \cite{LB}, we present an alternative strategy that can be applied to meromorphic functions.
In particular, working with the two transformations $\phi_1, \phi_2$ and indicating with $L$ the total number of points $(L=M+N+1)$, we show that the error decays like
\begin{equation*}
    \exp \left( - c \frac{L}{\left( \ln L \right)^k} \right), \quad L \rightarrow  + \infty,
\end{equation*}
with $k=1$ for $\phi_2$ ank $k=2$ for $\phi_1$, and where $c>0$ is a suitable constant.
Finally, we also propose a simple algorithm for automatic integration, that can be employed without information on the properties of the function $f$.

Throughout the work the symbol $\sim$ denotes the asymptotic equality, $\approx$ a generic approximation and $\lesssim$ states for less than or asymptotically equal to.

The paper is organized as follows.
In Section \ref{section 2} we give some general results regarding the trapezoidal rule and show how the error can be estimated in terms of contour integration.
In Section \ref{section 3} we recall the basic properties of the transformations $\phi_1$ and $\phi_2$.
Section \ref{section 4} deals with the error analysis of the trapezoidal rule based on the two transformations.
In Section \ref{section 5} we design an automatic integrator for Fourier type integrals and present some numerical examples.
Concluding remarks can be found in Section \ref{section 6}.

\section{General results for the trapezoidal rule} \label{section 2}

In this section we recall some theoretical results concerning the trapezoidal approximation
\begin{equation} \label{trap rule}
    I(f) = \int_{-\infty}^{+\infty} f(x) dx \approx h \sum_{j=-\infty}^{+\infty} f(jh),
\end{equation}
in which $f \colon \mathbb{R} \rightarrow \mathbb{R}$ is a generic integrable function and $h$ is a positive scalar.
Given $M$ and $N$ positive integers, we denote the truncated trapezoidal rule by
\begin{equation}
    T_{M,N,h}(f) = h \sum_{j=-M}^{N} f(jh).
\end{equation}
Then, for the quadrature error
\begin{equation*}
    \mathcal{E}_{M,N,h} = \left\vert I(f) - T_{M,N,h}(f) \right\vert,
\end{equation*}
it holds
\begin{equation*}
    \mathcal{E}_{M,N,h} \leq \mathcal{E}_D +\mathcal{E}_{T_L}+\mathcal{E}_{T_R},
\end{equation*}
where
\begin{align*}
    \mathcal{E}_D &= \left\vert \int_{-\infty}^{+\infty} f(x) dx - h \sum_{j=-\infty}^{+\infty} f(jh) \right\vert, \\
    \mathcal{E}_{T_L} &= h \left\vert \sum_{j=-\infty}^{-M-1}  f(jh) \right\vert, \quad
    \mathcal{E}_{T_R} = h \left\vert \sum_{j=N+1}^{+\infty}  f(jh) \right\vert.
\end{align*}
The quantities $\mathcal{E}_D$ and $\mathcal{E}_T:= \mathcal{E}_{T_L}+\mathcal{E}_{T_R}$ are referred to as the discretization error and the truncation error, respectively. We omit their dependence on $M,N,h$ to avoid redundant notations.

\begin{definition} \cite[Definition 2.12]{LB} Given $d>0$, let $\mathcal{D}_d$ be the infinite strip domain of width $2d$ given by
\begin{equation*}
   \mathcal{D}_d = \lbrace \zeta \in \mathbb{C} \colon | \Im (\zeta) | < d  \rbrace,
\end{equation*}
and let $\mathbf{B}\left(\mathcal{D}_d\right)$ be the set of functions analytic in $\mathcal{D}_d$ that satisfy
\begin{equation*}
    \int_{-d}^{d} | f(x+ i \eta ) | d \eta = \mathcal{O}\left( |x|^a \right), \quad x \rightarrow \pm \infty, \; 0 \leq a < 1,
\end{equation*}
and
\begin{equation*}
    \mathcal{N}(f,d)= \lim_{\eta \rightarrow d^-} \left\lbrace \int_{-\infty}^{+\infty} | f(x+ i \eta ) | d x + \int_{-\infty}^{+\infty} | f(x- i \eta ) | d x \right\rbrace < + \infty.
\end{equation*}
\end{definition}
For the discretization error of the trapezoidal rule applied to functions in $\mathbf{B}\left(\mathcal{D}_d\right)$, the following theorem holds (see \cite[Theorem 2.20]{LB}) .
\begin{theorem} 
    Assume $f \in \mathbf{B}\left(\mathcal{D}_d\right)$. Then 
    \begin{equation} \label{E_D old}
        \mathcal{E}_D \leq \frac{\mathcal{N}(f,d)}{2 \sinh(\pi d / h)} e^{-\frac{\pi d}{h}}.
    \end{equation}
\end{theorem}
By (\ref{E_D old}) we have that
\begin{equation*}
    \mathcal{E}_D \lesssim \mathcal{N}(f,d) e^{-\frac{2 \pi d}{h}}, \quad h \rightarrow 0.
\end{equation*}
The above result expresses the exponential decay of the discretization error as $h \rightarrow 0$, with a speed that increases with $d$. 
We remark that, if $f$ has a pole on the set $\left\vert \Im (\zeta) \right\vert= \overline{d}$, we may have $\mathcal{N}(f,d) \rightarrow + \infty$, for $d \rightarrow \overline{d}$.
Therefore, to obtain an accurate estimate, one should optimize the bound with respect to $0 < d < \overline{d}$, that in general may represent a difficult task.

In this work we follow a different approach, based on the analysis given in \cite{B} and \cite{DE}, that simplifies formula (\ref{E_D old}) when working with meromorphic functions.
First, let us consider a general quadrature rule 
\begin{equation} \label{quadrature lambda}
    I(f) = \int_{-\infty}^{+\infty} f(x) dx \approx \sum_{j=- \infty}^{+\infty} \lambda_j f \left(x_j \right),
\end{equation}
with remainder 
\begin{equation} \label{remainder}
    \mathcal{R}(f) =  I(f) - \sum_{j=-\infty}^{+\infty} \lambda_j f \left(x_j \right).
\end{equation}
Now, let $\varphi(z)$ and $\psi(z)$ be two analytic functions such that:
\begin{enumerate}
    \item[(i)] $\varphi(z)$ is a single-valued in the finite complex plane, without singularities and with distinct real zeros $\lbrace x_j \rbrace_{j \in \mathbb{Z}}$; 
    \item[(ii)] $\psi(z)$ has no singularities in the plane cut along the real axis;
    \item[(iii)] for $x \in \mathbb{R}$,
    \begin{equation*}
        \psi(x-0i)-\psi(x+0i) = 2 \pi i \varphi(x),
    \end{equation*}
    where
    \begin{equation*}
        \psi(x \pm 0i) = \lim_{y \rightarrow 0^+} \psi(x \pm i y).
    \end{equation*}
\end{enumerate}
By setting in (\ref{quadrature lambda})
\begin{equation*}
    \lambda_j =- \frac{\psi(x_j)}{\varphi^{'}(x_j)},
\end{equation*}
the remainder (\ref{remainder}) is given by
\begin{equation} \label{R(f)}
    \mathcal{R}(f) = \frac{1}{2 \pi i } \int_{\mathcal{C}} \frac{\psi(z)}{\varphi(z)} f(z) dz,
\end{equation}
where the contour $\mathcal{C}$ contains the real axis, but no singularity of the function $f(z)$ lies on or within the contour.

In \cite{B} and \cite{DE} it has been shown that the trapezoidal rule can be recast in this framework.
For $k>0$, we consider the functions
\begin{equation*}
    \varphi(z) = -\sin \left( z \sqrt{k} \right)
\end{equation*}
and
\begin{equation*}
    \psi(z) = 
    \begin{cases}
        \pi e^{i z \sqrt{k}}, \quad \Im(z) \geq 0\\
        \pi e^{-i z \sqrt{k}}, \quad \Im(z) <0
    \end{cases}.
\end{equation*}
These functions satisfy properties (i)-(ii)-(iii) and we observe that, for $k \rightarrow + \infty$,
\begin{equation} \label{psi su phi}
    \frac{\psi(z)}{\varphi(z)} \sim
    \begin{cases}
         2 \pi i e^{2iz \sqrt{k}}, \quad &\Im(z) \geq 0 \\
         -2 \pi i e^{-2iz \sqrt{k}}, \quad &\Im(z)<0 
    \end{cases}.
\end{equation}
Moreover, the zeros $\lbrace x_j \rbrace_{j \in \mathbb{Z}} $ of $\varphi(z)$ and the weights $\lambda_j$ are given by
\begin{equation*}
    x_j = \frac{j \pi}{\sqrt{k}} \quad {\rm and} \quad \lambda_j = \frac{\pi}{\sqrt{k}}.
\end{equation*}
By inserting these values in (\ref{quadrature lambda}) and by defining 
\begin{equation*}
    h:= \frac{\pi}{\sqrt{k}},
\end{equation*}
we obtain
\begin{equation*}
    I(f)  \approx h \sum_{j =-\infty}^{+\infty} f(jh),
\end{equation*}
that is, the trapezoidal rule (cf. (\ref{trap rule})).
In order to give an estimate of the remainder, let 
\begin{equation*}
    \mathcal{L}_R^{\pm} = \left\lbrace z \in \mathbb{C} \; \vert \; \Im(z) = \pm R \right\rbrace,
\end{equation*}
and suppose that $f(z)$ has no singularities on or between $\mathcal{L}_R^{+}$ and $\mathcal{L}_R^{-}$, excepts for a pair of simple poles $z_0$, and its conjugate $\bar{z_0}$.
Without loss of generality, we assume $\Im(z_0)>0$.
By choosing in (\ref{R(f)})
\begin{equation*}
    \mathcal{C} = \mathcal{C}_1 \cup \mathcal{C}_2 \cup \mathcal{L}_R^{+} \cup \mathcal{L}_R^{-},
\end{equation*}
where $\mathcal{C}_1$ and $\mathcal{C}_2$ are two small circles surrounding $z_0$ and $\Bar{z_0}$, we have that
\begin{equation*}
    \mathcal{R}(f) = \frac{1}{2 \pi i} \left[ \int_{\mathcal{L}_R^+ \cup \mathcal{L}_R^-} \frac{\psi(z)}{\varphi(z)} f(z) dz +\int_{\mathcal{C}_1 \cup \mathcal{C}_2} \frac{\psi(z)}{\varphi(z)} f(z) dz\right].
\end{equation*}
Now, by using (\ref{psi su phi}) the contribution given by the first integral is bounded by
\begin{equation} \label{int 1}
    \left\vert \frac{1}{2 \pi i} \int_{\mathcal{L}_R^+ \cup \mathcal{L}_R^-} \frac{\psi(z)}{\varphi(z)} f(z) dz \right\vert \lesssim e^{-2 R \sqrt{k}} \mathcal{N}(f,R).
\end{equation}
As for the second integral, by the residue theorem we obtain
\begin{equation*}
    \frac{1}{2 \pi i} \int_{\mathcal{C}_1 \cup \mathcal{C}_2} \frac{\psi(z)}{\varphi(z)} f(z) dz = - \left[ {\rm Res} \left( \frac{\psi(z)}{\varphi(z)} f(z), z_0 \right) + {\rm Res} \left( \frac{\psi(z)}{\varphi(z)}f(z), \Bar{z_0} \right) \right],
\end{equation*}
where the symbol ${\rm Res}( \cdot, \cdot)$ denotes the residue.
The minus sign in the above formula is due to the fact that, for travelling $\mathcal{C}$ in counterclockwise direction, the two circles are actually run clockwise.

Provided that the restriction of $f$ to the real numbers is real-valued, it holds 
\begin{equation*}
    {\rm Res} \left( \frac{\psi(z)}{\varphi(z)}f(z), z_0 \right) = \overline{{\rm Res} \left( \frac{\psi(z)}{\varphi(z)}f(z), \bar{z_0} \right)} = \frac{\psi(z_0)}{\varphi(z_0)}{\rm Res} \left( f(z),z_0 \right).
\end{equation*}
Remembering that $\sqrt{k}= \frac{\pi}{h}$, by using (\ref{psi su phi}) we finally have
\begin{align} \label{int 2}
    \frac{1}{2 \pi i} \int_{\mathcal{C}_1 \cup \mathcal{C}_2} \frac{\psi(z)}{\varphi(z)} f(z) dz &= -2 \frac{\psi(z_0)}{\varphi(z_0)} \Re \left\lbrace {\rm Res} \left( f(z),z_0 \right)\right\rbrace \notag \\
    &\sim -4 \pi \Re \left\lbrace i \rho_0 e^{2iz_0 \frac{\pi}{h}}  \right\rbrace, \quad h \rightarrow 0, \notag \\  
    &= 4 \pi \Im \left\lbrace \rho_0 e^{2iz_0 \frac{\pi}{h}}  \right\rbrace,
\end{align}
where $\rho_0 = {\rm Res}(f(z),z_0)$, and therefore
\begin{equation} \label{10 bis}
    \left\vert \frac{1}{2 \pi i} \int_{\mathcal{C}_1 \cup \mathcal{C}_2} \frac{\psi(z)}{\varphi(z)} f(z) dz \right\vert \lesssim 4 \pi \left\vert \rho_0 \right\vert e^{-2 \Im (z_0) \frac{\pi}{h}}.
\end{equation}
Since $R> \Im(z_0)$, by comparing (\ref{int 1}) and (\ref{10 bis}), we can neglect the contribution of the first integral and therefore we estimate $\mathcal{E}_D$ by 
\begin{equation} \label{E_D}
    \mathcal{E}_D \lesssim 4 \pi \left\vert \rho_0 \right\vert e^{-2 d \frac{\pi}{h}},
\end{equation}
in which $d =  \Im(z_0)$.
With respect to the bound (\ref{E_D old}), the above formula shows an error constant that is independent of $d$, and a faster exponential decay (in (\ref{E_D old}) one has to take $d < \Im(z_0)$).

Independently of the formula used to estimate $\mathcal{E}_D$, in order to design a reliable error approximation, one has to impose that $\mathcal{E}_D, \mathcal{E}_{T_L}, \mathcal{E}_{T_R}$ have the same exponential decay, in order to determine $h$, $M$, $N$, and then use these values to obtain an estimate of the total error $\mathcal{E}_{M,N,h}$.

\section{Double exponential transformations} \label{section 3}

In order to keep the exposition as clear as possible, we focus on the computation of the cosine transform
\begin{equation} \label{integrale coseno}
    F^{(c)}(\omega) = \int_0^{+\infty} f(x) \cos(\omega x)  dx, 
\end{equation}
and, where necessary, we simply explain which are the modifications needed to extend the analysis to the sine case.
As mentioned in Introduction, in order to numerically evaluate the above integral, we consider the transformations
\begin{equation} \label{DE}
    x =  \frac{\tau}{\omega} \phi_i\left(t -\frac{\pi}{2 \tau} \right), \quad  i=1,2, 
\end{equation}
where 
\begin{align}
    \phi_1(\xi) &= \frac{\xi}{1-e^{- 2 \pi \sinh \xi}}, \label{phi1} \\
    \phi_2(\xi) &= \frac{\xi}{1-e^{-2 \xi-\alpha \left( 1-e^{-t} \right) - \beta \left( e^t-1 \right)}}, \quad 0<\alpha<\beta<1, \label{phi2}
\end{align}
and $\tau>0$ is a given parameter.
The functions $\phi_i(\xi)$ are such that: 
\begin{itemize}
    \item[(a)] for $\xi \rightarrow - \infty$, $\phi_i (\xi) \rightarrow 0$ double exponentially, 
    \item[(b)] for $\xi \rightarrow + \infty$, $\phi_i(\xi) \rightarrow \xi$ double exponentially.
\end{itemize}
By using (\ref{DE}), integral (\ref{integrale coseno}) becomes
\begin{equation} \label{integrale trasf}
    F^{(c)}(\omega) = \frac{\tau}{\omega} \int_{-\infty}^{+\infty} f \left( \frac{\tau}{\omega} \phi_i\left(t -\frac{\pi}{2 \tau} \right) \right) \cos \left( \tau \phi_i\left(t -\frac{\pi}{2 \tau} \right) \right) \phi_i^{'}\left(t -\frac{\pi}{2 \tau} \right) d t,
\end{equation}
and the trapezoidal rule with mesh size $h$ reads
\begin{equation*}
    F^{(c)}_{h}(\omega) = \frac{\tau}{\omega} h \sum_{j = -\infty}^{+\infty} f \left( \frac{\tau}{\omega} \phi_i\left(jh -\frac{\pi}{2 \tau} \right) \right) \cos \left( \tau \phi_i\left(jh -\frac{\pi}{2 \tau} \right) \right) \phi_i^{'}\left(jh -\frac{\pi}{2 \tau} \right).
\end{equation*}
By taking $\tau$ such that $ \tau h = \pi$, and assuming $\left\vert f(x) \right\vert \leq C$, by properties (a) and (b) both tails decay double exponentially.

In order to work with the sine transform, we just need to use the change of variable $x = \frac{\tau}{\omega}\phi_i(t)$ in place of (\ref{DE}).

\section{Error analysis} \label{section 4}

In this section, we analyze the error of the truncated trapezoidal rule 
\begin{equation} \label{truncated}
 F^{(c)}_{M,N,h}(\omega) = \frac{\tau}{\omega} h \sum_{j = -M}^{N} f \left( \frac{\tau}{\omega} \phi_i\left(jh -\frac{\pi}{2 \tau} \right) \right) \cos \left( \tau \phi_i\left(jh -\frac{\pi}{2 \tau} \right) \right) \phi_i^{'}\left(jh -\frac{\pi}{2 \tau} \right),
\end{equation}
applied to (\ref{integrale trasf}).
The aim is to suitably define $h,M,N$ in order to equalize the error contributions $\mathcal{E}_D, \mathcal{E}_{T_L}, \mathcal{E}_{T_R}$.
As stated before, we assume $\left\vert f(x) \right\vert \leq C$, $x \in [0, +\infty)$.

\subsection{The transformation $\phi_1(\xi)$} \label{section phi 1}

First of all, by (\ref{phi1}) we have that 
\begin{equation} \label{phi asintotico}
    \phi_1(\xi) \sim 
    \begin{cases}
        | \xi | e^{2 \pi \sinh \xi}, \quad &\xi \rightarrow -\infty \\
        \frac{1}{2 \pi}, \quad &\xi \rightarrow 0 \\
        \xi, \quad &\xi \rightarrow + \infty
    \end{cases}.
\end{equation}
As for its derivative
\begin{equation*}
    \phi_1^{'}(\xi) = \frac{1- \left( 1+ 2 \pi \xi \cosh \xi \right) e^{-2 \pi \sinh \xi}}{\left(1-e^{-2 \pi \sinh \xi} \right)^2},
\end{equation*}
it holds
\begin{equation} \label{phi primo}
    \phi_1^{'}(\xi) \sim 
    \begin{cases}
        2 \pi | \xi | \cosh \xi e^{2 \pi \sinh \xi}, \quad &\xi \rightarrow -\infty \\
        \frac{1}{2}, \quad &\xi \rightarrow 0 \\
        1, \quad &\xi \rightarrow + \infty
    \end{cases}.
\end{equation}
We start the analysis by studying the truncation error $\mathcal{E}_{T_R}$, that is
\begin{equation*}
    \mathcal{E}_{T_R} = \left\vert \frac{\tau}{\omega} h \sum_{j = N+1}^{+\infty} f \left( \frac{\tau}{\omega} \phi_1\left(jh -\frac{\pi}{2 \tau} \right) \right) \cos \left( \tau \phi_1\left(jh -\frac{\pi}{2 \tau} \right) \right)  \phi_1^{'}\left(jh -\frac{\pi}{2 \tau} \right) \right\vert.
\end{equation*}
By (\ref{phi primo}), we have that
\begin{equation} \label{A}
    \mathcal{E}_{T_R} \lesssim \frac{\tau C}{\omega} h \sum_{j = N+1}^{+\infty} \left\vert \cos \left( \tau \phi_1\left(jh -\frac{\pi}{2 \tau} \right) \right) \right\vert, \quad N \rightarrow + \infty. 
\end{equation}
At this point we notice that, for $\xi \rightarrow + \infty$,
\begin{equation*}
    \cos \left( \tau \phi_1(\xi) \right) =\cos (\tau \xi)
    -\sin (\tau \xi) \left( \tau \left( \phi_1(\xi)- \xi \right) \right)+\mathcal{O} \left( \left( \phi_1(\xi)-\xi \right)^2 \right).
\end{equation*}
Now, since 
\begin{equation*}
    \tau \left( \phi_1(\xi) - \xi \right) = \tau \frac{\xi e^{-2 \pi \sinh \xi}}{1-e^{-2\pi \sinh \xi}} \leq \tau \cosh \xi e^{-2 \pi \sinh \xi},
\end{equation*}
by taking $\xi = jh-\frac{\pi}{2 \tau}$ and defining $\tau = \frac{\pi}{h}$, we have
\begin{equation*}
    \left\vert \cos \left( \tau \phi_1\left(jh -\frac{\pi}{2 \tau} \right) \right) \right\vert \leq \tau \cosh (jh) e^{-2 \pi \sinh (jh)}.
\end{equation*}
By inserting this result in (\ref{A}) we obtain
\begin{equation*}
    \mathcal{E}_{T_R} \lesssim \frac{\tau^2}{\omega} C h   \sum_{j = N+1}^{+\infty} \cosh(jh) e^{-2\pi \sinh(jh)} ,
\end{equation*}
and, finally, 
\begin{equation} \label{E_T_R}
    \mathcal{E}_{T_R}  \lesssim \frac{\tau^2 C}{2 \pi \omega} e^{-2\pi \sinh(Nh)}.
\end{equation}
As for the truncation error $\mathcal{E}_{T_L}$, we have that
\begin{align*}
    \mathcal{E}_{T_L}  &= \left\vert \frac{\tau}{\omega} h \sum_{j=-\infty}^{-M-1} f \left( \frac{\tau}{\omega} \phi_1\left(jh -\frac{\pi}{2 \tau} \right) \right) \cos \left( \tau \phi_1\left(jh -\frac{\pi}{2 \tau} \right) \right) \phi_1^{'}\left(jh -\frac{\pi}{2 \tau} \right) \right\vert \\
    &\leq \frac{\tau}{\omega} h C \sum_{j=-\infty}^{-M-1} \left\vert \phi_1^{'}\left(jh -\frac{\pi}{2 \tau} \right) \right\vert \\
    & \lesssim \frac{\tau C}{\omega} \phi_1 \left( -M h \right), \quad {\rm since} \; \phi_1^{'}(\xi) >0, \; \forall \xi \in \mathbb{R}.
\end{align*}
By using (\ref{phi asintotico}), we finally obtain 
\begin{equation} \label{E_T_L}
     \mathcal{E}_{T_L} \lesssim \frac{\tau C}{\omega} M h e^{-2 \pi \sinh (Mh)}, \quad {\rm for} \; M \rightarrow + \infty.
\end{equation}
As already mentioned, the idea now is to define $h,M,N$ such that $\mathcal{E}_D, \mathcal{E}_{T_R}, \mathcal{E}_{T_L}$ have the same exponential behavior.
By comparing (\ref{E_T_R}) and (\ref{E_T_L}) we simply impose $M=N$, for any given $N$.
As for the choice of $h$, we define it by solving the equation (see (\ref{E_D}), (\ref{E_T_R}))
\begin{equation} \label{impo}
    \frac{2 \pi d}{h} = 2 \pi \sinh(Nh).
\end{equation}
In the above formula, $d$ is given by the modulus of the imaginary part of the pole of the function
\begin{equation*}
    g(t) = f \left( \frac{\tau}{\omega} \phi_1 \left( t-\frac{\pi}{2 \tau} \right) \right),
\end{equation*}
closest to the real axis (cf. (\ref{E_D}) and (\ref{integrale trasf})). 
By \cite{OM2}, we have the following simple estimate
\begin{equation} \label{D Jap}
    d \sim   \frac{\theta}{\ln\left( \frac{\pi}{h} \right)} , \quad h \rightarrow 0,
\end{equation}
in which $\theta := \left\vert \arg(z_0) \right\vert$, where, as in Section \ref{section 2}, $z_0$ and its conjugate are the poles of $f$ closest to the real axis.
By inserting the above approximation in (\ref{impo}) and using 
\begin{equation*}
    \sinh(x) \sim \frac{e^{x}}{2}, \quad x \rightarrow + \infty,
\end{equation*}
we obtain
\begin{equation*}
  \ln\left( \frac{h}{\pi} \right) e^{\ln\left( \frac{h}{\pi} \right)} \sim - \frac{2 \theta}{\pi} e^{-Nh} ,
\end{equation*}
from which we have that 
\begin{equation*}
    \ln\left( \frac{h}{\pi} \right) \sim W_{-1} \left(  - \frac{2 \theta}{\pi} e^{-Nh} \right),  \quad N \rightarrow + \infty,
\end{equation*}
where $W_{-1}$ is the $-1$ branch of the Lambert $W$-function, sometimes also denoted by $W_m$ (see \cite[4.13]{NIST}).
Now, by using (see \cite[n. 4.13.11]{NIST})
\begin{equation*}
    W_{-1} (x) =  \ln(-x)-\ln \left(-\ln (-x) \right) + \mathcal{O} \left( 1 \right), \quad x \rightarrow 0^-,
\end{equation*}
for large $N$ we obtain
\begin{align*}
    \ln \left( \frac{h}{\pi} \right) &\sim \ln \left( \frac{2 \theta}{\pi}e^{-Nh} \right) - \ln \left( -\ln \left( \frac{2 \theta}{\pi} e^{-Nh} \right) \right) \\
    &= \ln \frac{\frac{2 \theta}{\pi}e^{-Nh}}{Nh -\ln \frac{2 \theta}{\pi}} \\
    &\sim \ln \frac{\frac{2 \theta}{\pi}e^{-Nh}}{Nh},
\end{align*}
and, therefore, 
\begin{equation*}
    h \sim \frac{2 \theta e^{-Nh}}{Nh}.
\end{equation*}
The above formula, after some computations, leads to
\begin{equation*}
    e^{\frac{Nh}{2}} \frac{Nh}{2} \sim \sqrt{\frac{N \theta}{2}},
\end{equation*}
and we finally obtain
\begin{equation} \label{h phi1}
    h \sim h^{\star} := \frac{2}{N} W_0 \left( \sqrt{\frac{N \theta}{2}} \right), \quad N \rightarrow + \infty,
\end{equation}
where $W_0$, also denoted by $W_p$, is the principal branch of the Lambert $W$-function (see \cite[n. 4.13.11]{NIST}).
At this point, the total error $\mathcal{E}_{M,N,h} = \mathcal{E}_D + \mathcal{E}_T$ is estimated as (see (\ref{E_D}), (\ref{E_T_R}), (\ref{E_T_L})) 
\begin{equation} \label{24bis}
  \mathcal{E}_{N,N,h^{\star}} \lesssim \left( \frac{\pi C}{2  \omega \left( h^{\star} \right)^2}+\frac{\pi C}{\omega} N\right) e^{-2 \pi \sinh \left( N h^{\star} \right)} +4 \pi |\rho_0|  e^{-2 d \frac{\pi}{h^{\star}}}.
\end{equation}
From the above formula we have that by increasing the frequency $\omega$ the error reduces.
Moreover, $\omega$ does not appear in the exponential terms and hence it does not affect the rate of convergence.

Denoting by $L = 2 N+1$ the total number of function evaluations, in the following proposition we show the asymptotic decay of the error with respect to $L$.
\begin{proposition}
For $L \rightarrow + \infty$, it holds
\begin{equation*}
    \mathcal{E}_{N,N,h^{\star}} \lesssim const \left( \frac{L}{\ln L} \right)^2  e^{-\pi \theta \frac{L}{\left( \ln L \right)^2}}.
\end{equation*}
\end{proposition}
\begin{proof}
    Denoting by  
    \begin{equation*}
        p_D(h) = \frac{2 \pi d}{h} \quad {\rm and} \quad p_{T_R} (N,h) = 2 \pi \sinh (Nh),
    \end{equation*}
    (cf. (\ref{E_D})-(\ref{E_T_R})), by using (\ref{h phi1}) and the approximation (see \cite[n. 4.13.11]{NIST})
    \begin{equation} \label{W_0 aprox}
        W_0(x) = \ln x \left( 1 + \mathcal{O} \left( \frac{\ln (\ln x)}{\ln x} \right) \right), \quad x \rightarrow + \infty,
    \end{equation}
    by direct computation we find
    \begin{equation} \label{pD pTR}
        p_D \left( h^{\star} \right), p_{T_R}\left(N,h^{\star} \right) = 2 \pi \theta \frac{N}{\left( \ln N \right)^2} \left( 1 + \mathcal{O} \left( \frac{\ln (\ln N)}{\ln N} \right) \right).
    \end{equation}
    By inserting (\ref{pD pTR}) in (\ref{24bis}), using again (\ref{h phi1}) for the error constant and since $N = \frac{L-1}{2}$, we obtain the result.
\end{proof}

\begin{figure}
    \begin{center}
    \includegraphics[scale=0.35]{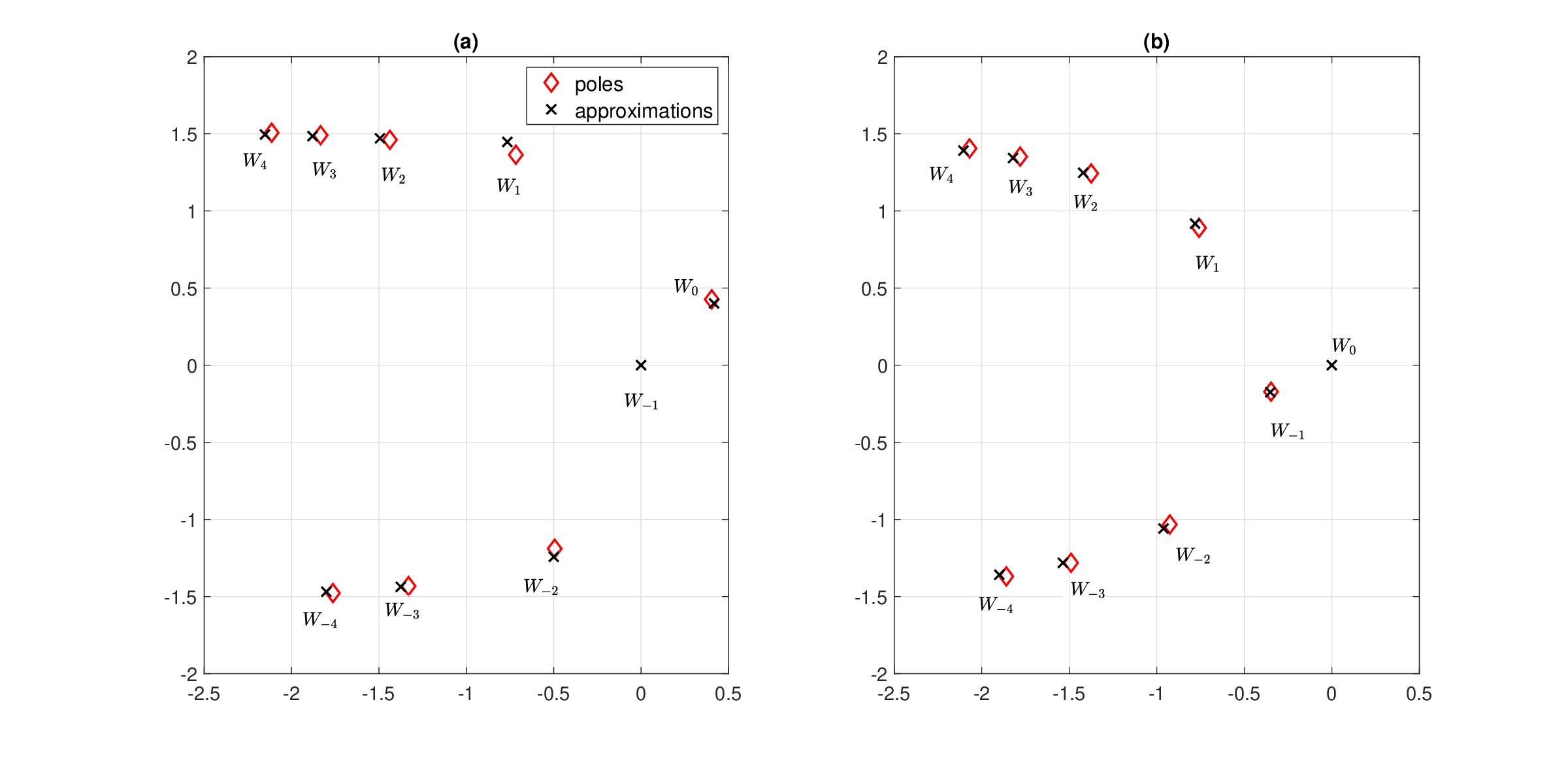}
    \end{center}
    \caption{The poles of the function $g(t)$ (red diamonds) and the approximations (\ref{tjk}) (black crosses) for $z_0 = 1+i$, $\omega = 4$, $h=0.3$ (left) and for $z_0 = 1-i$, $\omega = 2$, $h=0.05$ (right). In order to approximate $t^{(0)}$, the correct branch to use in (\ref{tjk}) is $W_0$ in (a) and $W_{-1}$ in (b).} \label{figura_poli}
\end{figure}

\begin{figure}
    \begin{center}
    \includegraphics[scale=0.35]{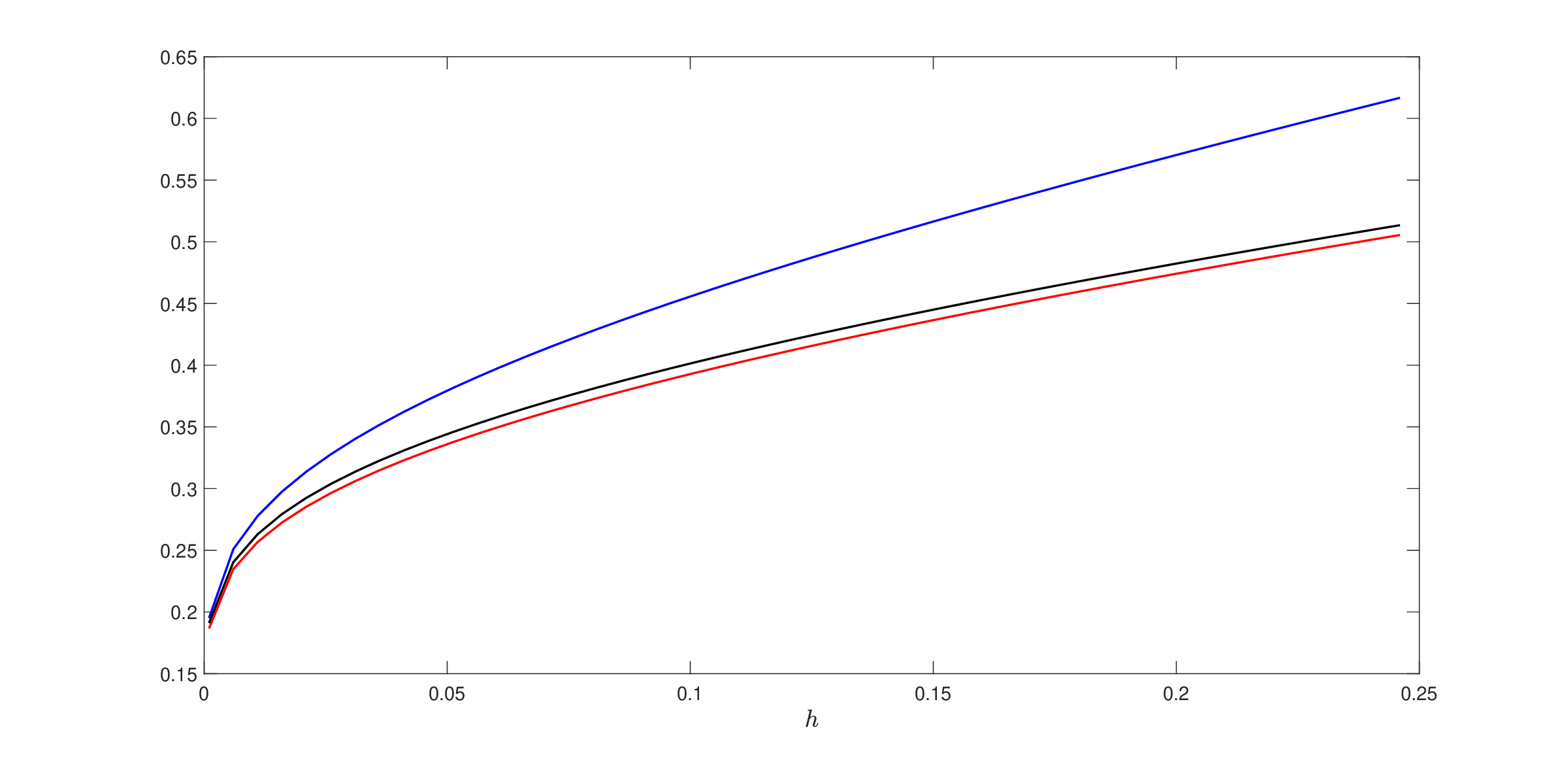}
    \end{center}
    \caption{The imaginary part of $t^{(0)}$ (red), the imaginary part of the approximation (\ref{tjk}) (black) and approximation (\ref{D Jap}) (blue) for different values of $h$. In this case $z_0 = i$ and $\omega=3$.} \label{figura_poli_H}
\end{figure}

The approximation (\ref{D Jap}) has been used to obtain a simple and fairly good approximation of the optimal value of $h$ for a given $N$ (cf. (\ref{h phi1})).
Anyway, in the numerical experiments to test the accuracy of (\ref{24bis}) we have used a different approximation of $d$.
In particular, by looking for the poles of
\begin{equation*}
    g(t) = f \left( \frac{\tau}{\omega} \phi_1 \left( t-\frac{\pi}{2 \tau} \right) \right), \quad \tau = \frac{\pi}{h},
\end{equation*}
we have to solve $g(t) = z_0$, that is, with respect to $\xi = t-\frac{h}{2}$,
\begin{equation*}
    \frac{\xi}{1-e^{-2 \pi \sinh \xi}} = z_0 \frac{\omega}{\pi} h,
\end{equation*}
(cf. (\ref{phi1})).
In this setting, in order to reduce the complexity of the above equation, we replace it by
\begin{equation*}
    \frac{\sinh \xi}{1-e^{-2 \pi \sinh \xi}} = z_0 \frac{\omega}{\pi} h,
\end{equation*}
to obtain
\begin{equation*}
    \sinh \left( \xi \right) = z_0 \frac{\omega}{\pi} h  +\frac{1}{2 \pi} W \left( 
-2  z_0 \omega h e^{-2 z_0\omega h} \right),
\end{equation*}
and therefore
\begin{equation} \label{tjk}
    \xi = {\rm arcsinh} \left( z_0 \frac{\omega}{\pi} h  +\frac{1}{2 \pi} W \left( 
-2 z_0 \omega h e^{-2 z_0 \omega h} \right) \right).
\end{equation}
Since both $\sinh ( \cdot)$ and $W( \cdot)$  are multifunctions, we have a set $\left\lbrace t_{jk} \right\rbrace_{j,k \in \mathbb{Z}}$ of poles.
The closest to the real axes, that we denote by $t^{(0)}$, is obtained by choosing the principal branch of $\sinh$ and, for $h$ sufficiently small, the branch $-1$ or $1$ of the Lambert $W$-function.
Without any assumption on $h$, the general rule for the choice of the branch of $W$ is the following. 
For $n \in \mathbb{Z}$, let $\mathcal{W}_n$ be the region of the complex plane defined by 
\begin{equation*}
    \mathcal{W}_n = \left\lbrace z \in \mathbb{C} \mid W_n \left( z e^z \right) = z \right\rbrace
\end{equation*}
(see \cite[Figure 4.13.2]{NIST}).
Defining $q = -2z_0\omega h$ (cf. (\ref{tjk})), if $\Im \left( z_0 \right)>0$, then $q \in \mathcal{W}_n$ for a certain $n \leq 0$ and we have to consider the branch $W_{n+1}$.
On the other side, if $\Im \left( z_0 \right)<0$, then $q \in \mathcal{W}_n$ for a suitable $n \geq 0$ and we have to take the branch $W_{n-1}$.
In both cases, $n = 0$ as $h \rightarrow 0$ and the branches to consider are $-1$ or $1$, depending on the sign of $\Im (z_0)$.
We observe that whenever
\begin{equation*}
    W_n\left(qe^q \right) = q,
\end{equation*}
we obtain $\xi = {\rm arcsinh} (0) = 0$, that is not an acceptable solution.
In Figure \ref{figura_poli} two examples are reported.
In Figure \ref{figura_poli_H} we have also compared the imaginary part of $t^{(0)}$ (numerically computed) together with the imaginary part of the approximations (\ref{tjk}) and (\ref{D Jap}).


\subsection{The transformation $\phi_2(\xi)$}

As pointed out in \cite{OM2}, the basic drawback of the transformation $\phi_1$ is that $d \rightarrow 0$ for $h \rightarrow 0$ (cf. (\ref{D Jap})) and this clearly slows down the method. 
This issue has motivated the introduction of the transformation 
\begin{equation*}
    \phi_2(\xi) = \frac{\xi}{1-e^{-2 \xi-\alpha \left( 1-e^{-\xi} \right) - \beta \left( e^{\xi}-1 \right)}},
\end{equation*}
with the suggested values
\begin{equation*}
    \beta = \frac{1}{4}, \quad \alpha = \frac{\beta}{\sqrt{1+\frac{\tau}{4 \pi} \log \left( 1+ \tau \right)}}.
\end{equation*}
We first observe that 
\begin{equation} \label{phi_2 asympt}
    \phi_2(\xi) \sim 
    \begin{cases}
        |\xi| e^{-\alpha e^{-\xi}}, \quad &\xi \rightarrow -\infty,\\
        \frac{1}{2+\alpha +\beta}, \quad &\xi \rightarrow 0, \\
        \xi, \quad &\xi \rightarrow + \infty.
    \end{cases}
\end{equation}
Moreover, for $h \rightarrow 0$,
\begin{align} \label{alpha h piccolo}
    \alpha &= \frac{\beta}{\sqrt{1+\frac{\tau}{4 \pi}\ln (1+\tau)}} = \frac{\beta}{\sqrt{1+\frac{1}{4 h}\ln \left(1+\frac{\pi}{h}\right)}} \notag \\
    &= 2 \beta \sqrt{\frac{h}{-\ln h}} \left( 1+\mathcal{O}\left( \frac{1}{\ln h} \right) \right).
\end{align}
As before, we start by studying the truncation error $\mathcal{E}_{T_R}$. 
By exploiting the previous analysis (see (\ref{A})) and since, for $\xi > 0$ sufficiently large,
\begin{equation*}
    \tau (\phi_2(\xi)-\xi) = \tau \frac{\xi e^{-2\xi-\alpha(1-e^{-\xi})-\beta(e^{\xi}-1)}}{1-e^{-2\xi-\alpha(1-e^{-\xi})-\beta(e^{\xi}-1)}}\leq \tau \cosh \xi e^{-2 \beta \sinh \xi},
\end{equation*}
we obtain 
\begin{equation} \label{E_T_L_2}
    \mathcal{E}_{T_R} \lesssim \frac{\tau^2 C}{2 \beta \omega} e^{-2\beta \sinh (Nh)}, \quad N \rightarrow + \infty.
\end{equation}
As for the truncation error $\mathcal{E}_{T_L}$, by using (\ref{phi_2 asympt}), and following the analysis given for $\phi_1$, we obtain
\begin{equation} \label{E_T_R_2}
     \mathcal{E}_{T_L} \lesssim \frac{\tau C}{\omega}  M h e^{-\alpha e^{Mh}}, \quad M \rightarrow + \infty.
\end{equation}
At this point, in order to determine $h$, we consider $\mathcal{E}_D$ and $\mathcal{E}_{T_R}$ as in (\ref{E_D}) and (\ref{E_T_L_2}), respectively, and impose
\begin{equation*} 
    2 \beta \sinh (Nh) = 2 d \frac{\pi}{h}.
\end{equation*}
As for the value of $d$, we use the estimate given in \cite{OM2}, that is,
\begin{equation*}
    d \sim \frac{\theta}{2}, \quad h \rightarrow 0,
\end{equation*}
$\theta = \left\vert \arg (z_0) \right\vert$, and therefore
\begin{equation} \label{h segnato}
    h \sim \overline{h} := \frac{1}{N} W_0 \left( \frac{N  \pi \theta}{\beta} \right), \quad N \rightarrow + \infty.
\end{equation}
By considering the truncation errors $\mathcal{E}_{T_R}$ and $\mathcal{E}_{T_L}$ (see (\ref{E_T_L_2}), (\ref{E_T_R_2})) and imposing 
\begin{equation*}
    2 \beta \sinh (Nh) = \alpha e^{Mh},
\end{equation*}
for any given $N$ we define
\begin{equation} \label{M ceil}
    M = \left\lceil N-\frac{1}{\overline{h}} \ln \left( \frac{\alpha}{\beta} \right) \right\rceil ,
\end{equation}
as approximate solution, where $\left\lceil \; \right\rceil$ denotes the ceil operator.
Note that $M \geq N$, since $\beta > \alpha$.
The total error is finally estimated by
\begin{equation} \label{errore phi2}
    \mathcal{E}_{M,N,\overline{h}} \lesssim \frac{\pi^2 C}{2 \beta \omega \left( \overline{h}\right)^2} e^{-2 \beta \sinh \left( N \overline{h} \right)} +\frac{\pi C}{\omega}  M e^{-\alpha e^{M \overline{h}}} +4 \pi |\rho_0|   e^{-2 d \frac{\pi}{\overline{h}}}.
\end{equation}
As before (cf. formula (\ref{24bis})), $\omega$ does not appear in the exponential terms and the error reduces for growing $\omega$.

Denoting by $L=M+N+1$ the total number of points, we have the following result.
\begin{proposition}
    For $L \rightarrow + \infty$, it holds
    \begin{equation*}
        \mathcal{E}_{M,N,\overline{h}} \lesssim const \left( \frac{L}{\ln L} \right)^2 e^{-\frac{2}{5} \theta \pi \frac{L}{\ln L}}.
    \end{equation*}
\end{proposition}
\begin{proof}
    Denoting by 
    \begin{equation*}
        q_D(h) = \frac{2 \pi d}{h}, \;\; q_{T_R}(N,h) = 2 \beta \sinh (Nh), \;\; q_{T_L}(M,h) = \alpha e^{Mh},
    \end{equation*}
    (cf. (\ref{E_D})-(\ref{E_T_L_2})-(\ref{E_T_R_2})), using (\ref{h segnato}), (\ref{M ceil}) and (\ref{W_0 aprox}), by direct computation we find
    \begin{equation} \label{q q q}
        q_D(\overline{h}), q_{T_R}(N,\overline{h}), q_{T_L}(M,\overline{h}) = \theta \pi \frac{N}{\ln N} \left( 1 + \mathcal{O} \left( \frac{\ln (\ln N)}{\ln N} \right) \right).
    \end{equation}
    Now, by definition (\ref{M ceil}) and inserting (\ref{h segnato}) in (\ref{alpha h piccolo}), for $N \rightarrow + \infty$ we obtain 
    \begin{align*}
        M &\sim N -\frac{N}{W_0 \left( \frac{ \theta \pi}{\beta} N \right)} \ln \left( 2 \sqrt{\frac{W_0 \left( \frac{ \theta \pi}{\beta} N \right)}{-N \ln \frac{W_0 \left( \frac{ \theta \pi}{\beta} N \right)}{N}}} \right) \\
        &= N -\frac{N}{W_0 \left( \frac{ \theta \pi}{\beta} N \right)} \ln \left( 2 \sqrt{\frac{W_0 \left( \frac{ \theta \pi}{\beta} N \right)}{-N  \left( -W_0 \left( \frac{ \theta \pi}{\beta}  \right) +\ln \left( \frac{ \theta \pi}{\beta} \right)\right)}} \right) \\
        &\sim N -\frac{N}{W_0 \left( \frac{ \theta \pi}{\beta} N \right)} \left( -\frac{W_0\left( \frac{ \theta \pi}{\beta} N \right)}{2}+ \frac{1}{2} \ln \frac{ \frac{ \theta \pi}{\beta}}{W_0\left( \frac{ \theta \pi}{\beta} N \right)} \right) \\ 
        &\sim  N - \frac{N}{W_0 \left( \frac{ \theta \pi}{\beta} N \right)} \left( -\frac{1}{2} \ln N \right) \\
        &\sim N+\frac{N}{2} \frac{\ln N}{\ln \frac{N \theta \pi}{\beta}} \\
        &\sim \frac{3}{2}N,
    \end{align*}
    where we have also used the relation (see \cite[n.4.13.5]{NIST})
    \begin{equation*}
        \ln \frac{x}{W_0(x)} = W_0(x),
    \end{equation*}
    and (\ref{W_0 aprox}).
    Then, for $L=N+M+1$ we have
    \begin{equation} \label{N caso peggiore}
        L \lesssim \frac{5}{2} N, \quad {\rm for} \quad N \rightarrow + \infty.
    \end{equation}    
    By inserting (\ref{N caso peggiore}) and (\ref{q q q}) in (\ref{errore phi2}) and using again (\ref{h segnato}) for the error constant, we obtain the result.
\end{proof}

\subsection{Numerical examples}

\begin{figure}
    \begin{center}
    \includegraphics[scale=0.35]{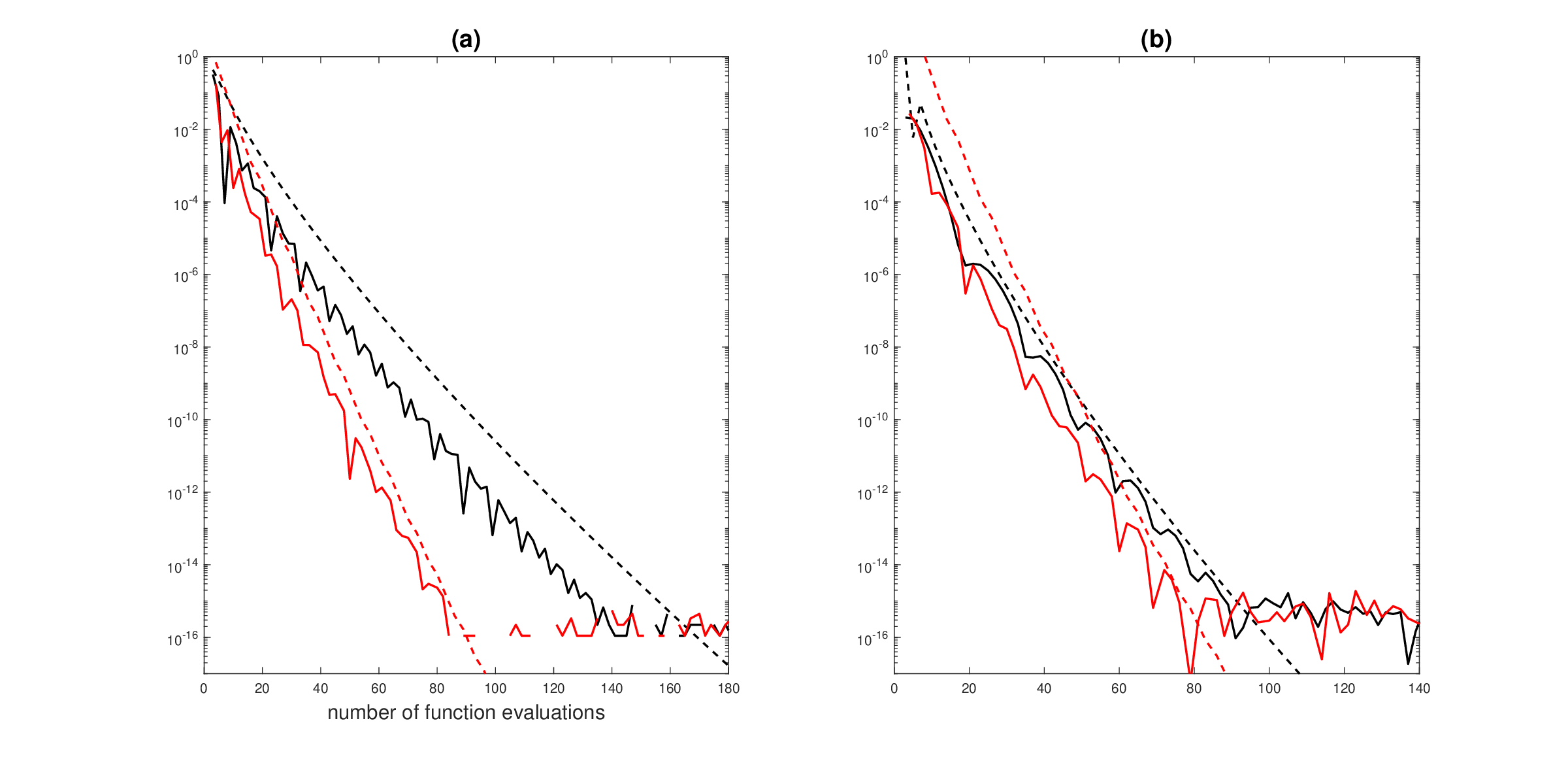}
    \end{center}
    \caption{The absolute error of the trapezoidal rule based on the transformations $\phi_1$ (black line) and $\phi_2$ (red line) with the estimates (\ref{24bis}) (dashed black line) and (\ref{errore phi2}) (dashed red line). The function considered is $f_1(x)$, with $\omega=1$ (left) and $\omega=7$ (right).} \label{figura_errori1}
\end{figure}

\begin{figure}
    \begin{center}
    \includegraphics[scale=0.35]{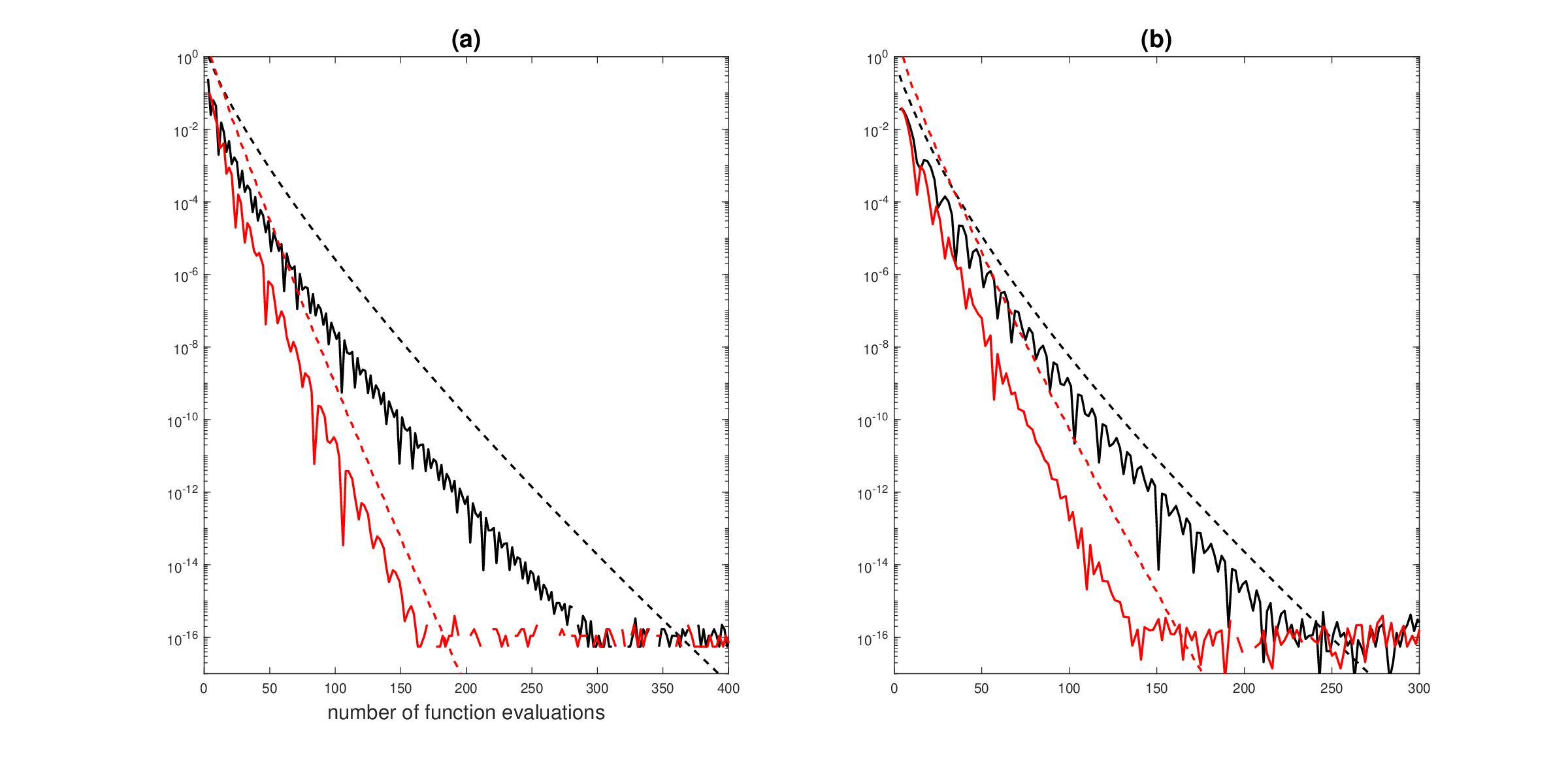}
    \end{center}
    \caption{The absolute error of the trapezoidal rule based on the transformations $\phi_1$ (black line) and $\phi_2$ (red line) with the estimates (\ref{24bis}) (dashed black line) and (\ref{errore phi2}) (dashed red line). The function considered is $f_2(x)$, with $\omega=1.5$ (left) and $\omega=5$ (right).} \label{figura_errori2}
\end{figure}

\begin{figure}
    \begin{center}
    \includegraphics[scale=0.35]{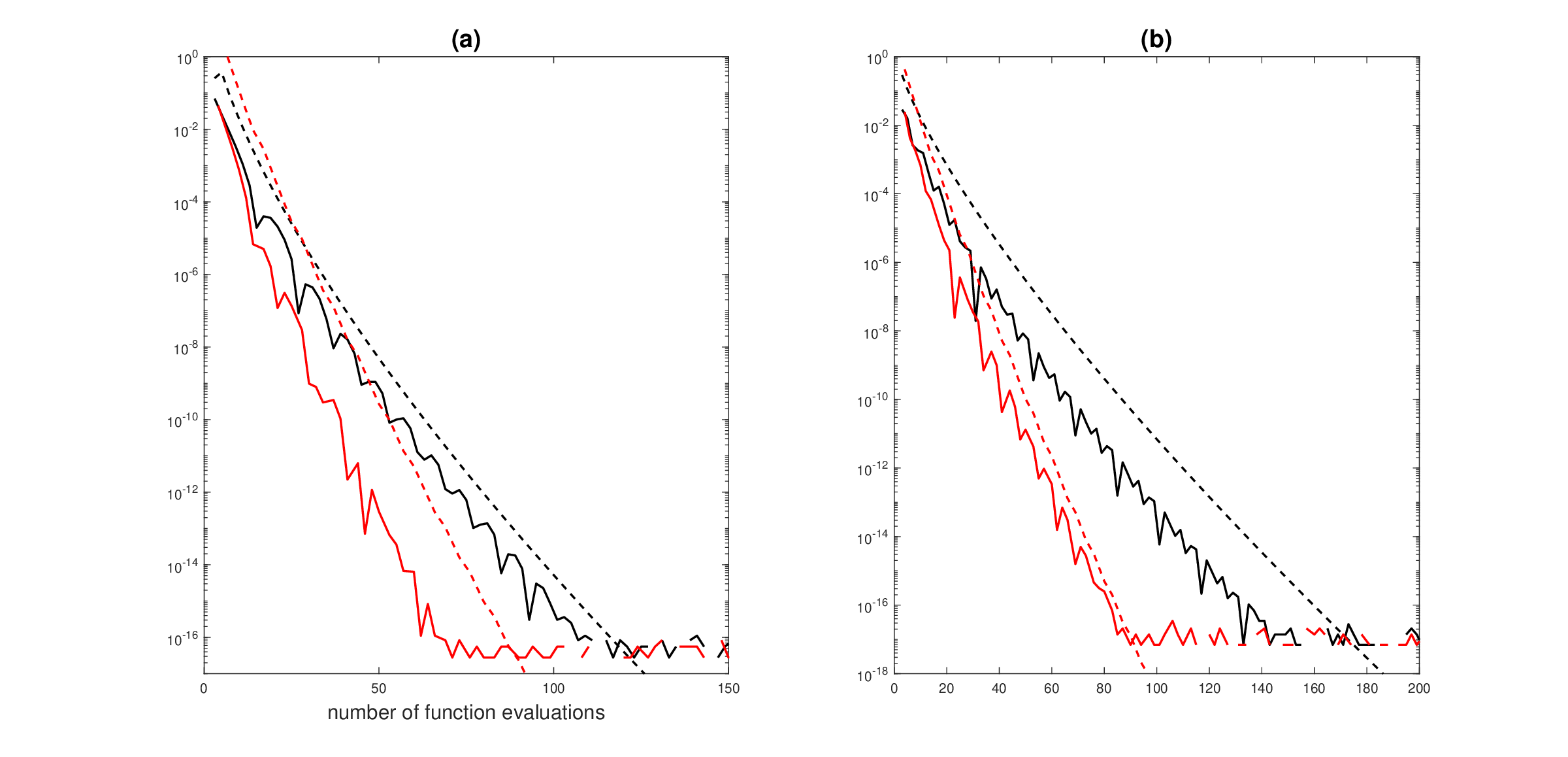}
    \end{center}
    \caption{The absolute error of the trapezoidal rule based on the transformations $\phi_1$ (black line) and $\phi_2$ (red line) with the estimates (\ref{24bis}) (dashed black line) and (\ref{errore phi2}) (dashed red line). The function considered is $f_3(x)$, with $\delta=1.5$ (left) and $\delta=5$ (right). In both cases $\omega = 2$.} \label{figura_errori3}
\end{figure}

\begin{figure}
    \begin{center}
    \includegraphics[scale=0.35]{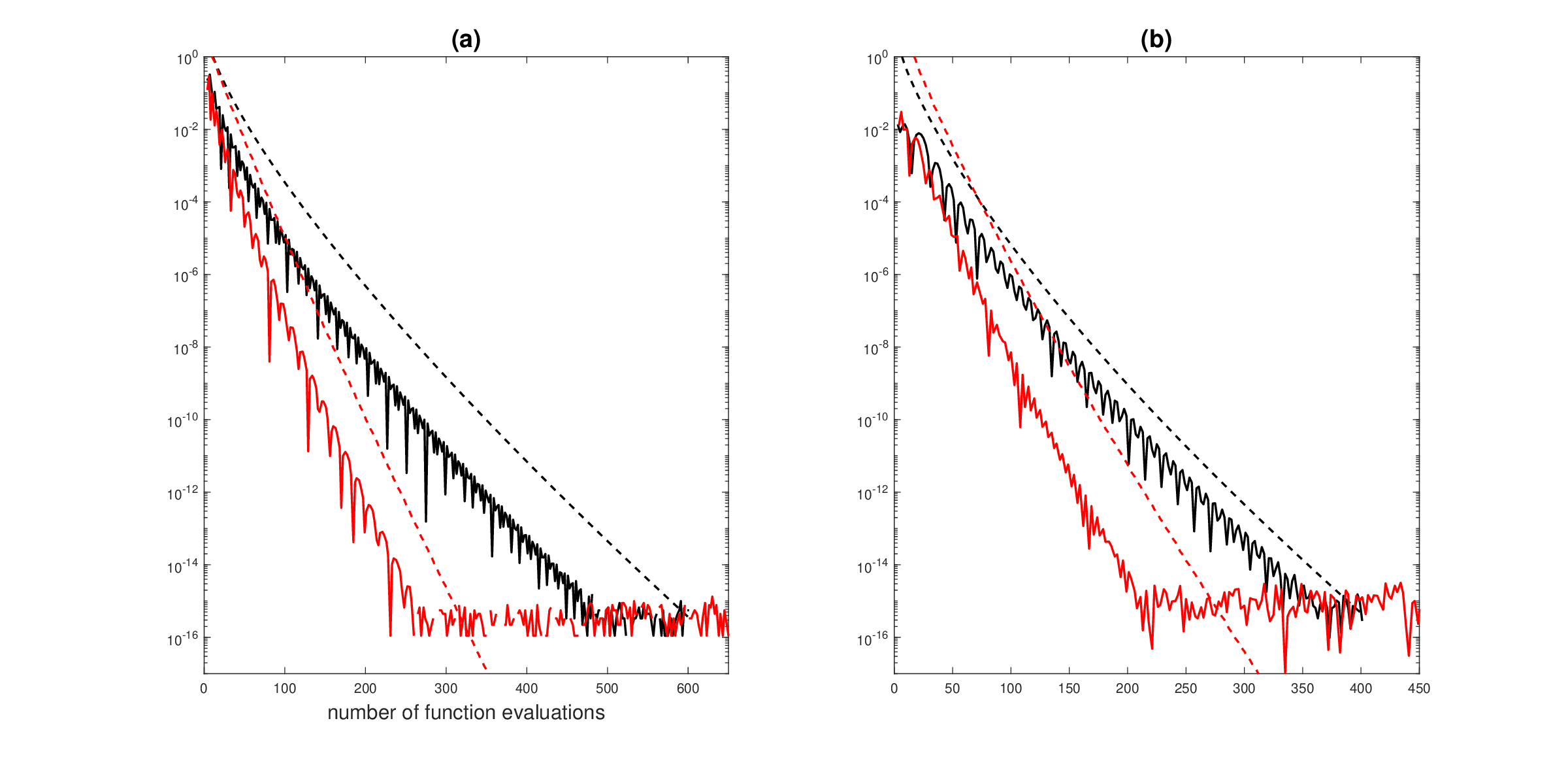}
    \end{center}
    \caption{The absolute error of the trapezoidal rule based on the transformations $\phi_1$ (black line) and $\phi_2$ (red line) with the estimates (\ref{24bis}) (dashed black line) and (\ref{errore phi2}) (dashed red line). The function considered is $f_4(x)$, with $\omega=1$ (left) and $\omega=4$ (right).} \label{figura_errori4}
\end{figure}

In this section we consider some numerical experiments in which we compare the errors of the trapezoidal rules based on the two transformations, together with the error estimates derived in the previous sections.
The functions considered are 
\begin{align*}
    f_1(x) &= \frac{1}{1+x^2}, \quad {\rm with \; poles} \quad z_k = (-1)^k i, \; k=0,1,\\
    f_2(x) &= \frac{x}{1+x^4}, \quad {\rm with \; poles} \quad z_k = e^{i \frac{\pi}{4} \left( 1+2k \right)}, \; k=0,1,2,3,\\
    f_3(x) &= \frac{1}{1+e^{\delta x}}, \quad {\rm with \; poles} \quad z_k = i(2k+1) \frac{\pi}{\delta}, \; k \in \mathbb{Z},\\
    f_4(x) &= \frac{1}{(x-2)^2+1}, \quad {\rm with \; poles} \quad z_k = 2+(-1)^k i, \; k=0,1.
\end{align*}
In Figures \ref{figura_errori1}-\ref{figura_errori2}-\ref{figura_errori3}-\ref{figura_errori4}, for different values of $\omega$, we plot the absolute error of the trapezoidal rule based on the transformations $\phi_1$ and $\phi_2$, and the error estimates given by formulas (\ref{24bis}) and (\ref{errore phi2}).
We can see that, especially for low frequencies, the trapezoidal rule based on the transformation $\phi_2$ provides better results, while for higher frequencies the difference between the two rules is less evident (see Figure \ref{figura_errori1}).
Anyway, the plots show the good accuracy of both error estimates.

\section{Automatic integration} \label{section 5}

In this section we design an automatic integrator for the computation of the sine and cosine transforms $F^{(s)}(\omega)$ and $F^{(c)}(\omega)$, by using the transformation $\phi_1$ of Section \ref{section phi 1}.
The idea is the following.
Setting an arbitrary tolerance $\eta$, we first look for $\ell$ such that the contribute of the truncation error $\mathcal{E}_T$ is less than $\frac{2}{3}\eta$ in $\mathbb{R} \setminus \left[ -\ell, \ell \right]$.
Exploiting the fact that, taking $M=N$, $\mathcal{E}_{T_R}$ and $\mathcal{E}_{T_L}$ have the same exponential decay, we consider the estimate (see formula (\ref{E_T_R}))
\begin{equation} \label{epsilon}
    \mathcal{E}_{T_R} \lesssim \frac{\tau^2 C}{2 \pi \omega} e^{-2 \pi \sinh(Nh)}.
\end{equation}
Then, we define $\ell = Nh$ by solving 
\begin{equation*}
    e^{-2 \pi \sinh \ell} = \frac{\eta}{3},
\end{equation*}
that is,
\begin{equation} \label{L}
    \ell ={\rm arcsinh} \left( -\frac{1}{2 \pi} \ln \frac{\eta}{3} \right).
\end{equation}
In this way, the truncation error is approximately given by 
\begin{equation*}
    \mathcal{E}_T = \mathcal{E}_{T_R} + \mathcal{E}_{T_L} \approx \frac{2}{3} \eta. 
\end{equation*}
Now, in order to properly define $N$ and $h$, we try to optimize the choice of $d$ in the discretization error $\mathcal{E}_D$, given by formula (\ref{E_D old}), as follows.
Taking $\Tilde{\ell} = \ell \gamma$, $\gamma>1$, we consider the estimate 
\begin{equation} \label{E_D delta}
    \mathcal{E}_D \approx \left\vert F^{(\cdot)}_{N_1,N_1,h_1} - F^{(\cdot)}_{N_2,N_2,h_2} \right\vert =: \Delta,
\end{equation}
where $N_2 = 2 N_1$ (in order to work with an embedded formula for computing $\Delta$), $h_1 = \frac{\Tilde{\ell}}{N_1}$ and $h_2 = \frac{h_1}{2}$.
The use of $\Tilde{\ell} >\ell$ is to make dominant the discretization error and therefore more reliable the approximation (\ref{E_D delta}).
In our examples we have set $\gamma=1.2$ and $N_1 =10 \div 20$.
Since the discretization error decays like
\begin{equation} \label{D}
   \mathcal{E}_D \sim c e^{-2 \pi \frac{d}{h}},
\end{equation}
we define $d$ by imposing
\begin{equation*}
    e^{- 2 \pi \frac{d}{h_1}} = \Delta,
\end{equation*}
that is,
\begin{equation} \label{d automatico}
    d = -\frac{h_1}{2 \pi} \ln \Delta.
\end{equation}
By comparing the exponential terms of (\ref{epsilon}) and (\ref{D}), we find 
\begin{equation*}
   h \sim h_d := \frac{1}{N} W_0 \left( 2 N d\right), \quad N \rightarrow + \infty. 
\end{equation*}
Since for our purposes $N h = \ell$, by imposing $h=h_d$, we obtain $W_0(2Nd)=\ell$, that leads to
\begin{equation*}
    \ell e^{\ell}=2Nd.
\end{equation*}
Therefore, we define
\begin{equation*}
    N = \left\lceil \frac{\ell e^{\ell}}{2 d} \right\rceil.
\end{equation*}
We summarize the above strategy in the following algorithm.
\begin{algorithm} \label{algoritmo}
    Given $\eta >0$, $N_1$, $N_2=2N_1$, $\gamma>1$,
    \begin{enumerate}
        \item compute $\ell$ as in (\ref{L}),
        \item compute $d$ as in (\ref{d automatico}), with 		$\Delta$ as in (\ref{E_D delta}),
        \item set $N = \left\lceil \frac{\ell e^{\ell}}{2 d} \right\rceil$ and $M=N$,
        \item set $h = \frac{\ell}{N}$.
    \end{enumerate}
\end{algorithm}
In Tables \ref{table1}-\ref{table2}-\ref{table3}-\ref{table4}, working with $\eta = 1e-7$, $\eta = 1e-10$ and $\eta = 1e-13$, and different values of $\omega$, we test Algorithm \ref{algoritmo} on some examples.
In particular, we report the final error, obtained by considering a reference solution, together with the corresponding values of $h$ and $N$.

\begin{table} 
\[
\begin{array}{cccccccccccc}
\toprule
\multicolumn{2}{c}{} & \multicolumn{1}{c}{\eta=1e-7} & \multicolumn{1}{c}{\eta=1e-10} & \multicolumn{1}{c}{\eta=1e-13} \\
\midrule
 & {\rm error} & 2.78e-7 & 1.68e-10 & 2.19e-13 \\
\omega = 1 & h & 9.63e-2 & 7.61e-2 & 5.75e-2 \\
& N & 18 & 27 & 40 \\
\midrule
 & {\rm error} & 4.66e-7 & 2.84e-11 & 1.01e-13 \\
\omega = 5 & h & 1.58e-1 & 9.34e-2 & 6.97e-2 \\
& N & 11 & 22 & 33 \\
\midrule
 & {\rm error} & 6.81e-9 & 2.39e-11 & 1.55e-14 \\
\omega = 10 & h & 1.24e-1 & 8.93e-2 & 7.42e-2 \\
& N & 14 & 23 & 31 \\
\bottomrule
\end{array}
\]
\caption{Results of Algorithm \ref{algoritmo} for the cosine transform $F^{(c)}(\omega)$, with $f(x) = \frac{1}{1+x^2}$ and $N_1=10$.} \label{table1}
\end{table}

\begin{table} 
\[
\begin{array}{cccccccccccc}
\toprule
\multicolumn{2}{c}{} & \multicolumn{1}{c}{\eta=1e-7} & \multicolumn{1}{c}{\eta=1e-10} & \multicolumn{1}{c}{\eta=1e-13} \\
\midrule
 & {\rm error} & 6.19e-8 & 1.05e-10 & 2.03e-12 \\
\omega = 1 & h & 5.25e-2 & 3.67e-2 & 2.91e-2 \\
& N & 33 & 56 & 79 \\
\midrule
 & {\rm error} & 3.29e-8 & 1.33e-10 & 5.14e-12 \\
\omega = 5 & h & 6.19e-2 & 4.67e-2 & 3.90e-2 \\
& N & 28 & 44 & 59 \\
\midrule
 & {\rm error} & 1.19e-8 & 2.94e-10 & 9.82e-16 \\
\omega = 10 & h & 7.22e-2 & 5.71e-2 & 4.26e-2 \\
& N & 24 & 36 & 54 \\
\bottomrule
\end{array}
\]
\caption{Results of Algorithm \ref{algoritmo} for the sine transform $F^{(s)}(\omega)$, with $f(x) = \frac{x}{1+x^4}$ and $N_1=20$.} \label{table2}
\end{table}

\begin{table} 
\[
\begin{array}{cccccccccccc}
\toprule
\multicolumn{2}{c}{} & \multicolumn{1}{c}{\eta=1e-7} & \multicolumn{1}{c}{\eta=1e-10} & \multicolumn{1}{c}{\eta=1e-13} \\
\midrule
 & {\rm error} & 9.17e-8 & 4.71e-11 & 5.54e-14 \\
\omega = 1 & h & 1.24e-1 & 7.90e-2 & 6.77e-2 \\
& N & 14 & 26 & 34 \\
\midrule
 & {\rm error} & 1.61e-9 & 7.06e-12 & 1.08e-14 \\
\omega = 5 & h & 1.24e-1 & 9.34e-2 & 7.42e-2 \\
& N & 14 & 22 & 31 \\
\midrule
 & {\rm error} & 1.43e-8 & 3.75e-9 & 2.39e-13 \\
\omega = 10 & h & 1.33e-1 & 1.21e-1 & 8.52e-2 \\
& N & 13 & 17 & 27 \\
\bottomrule
\end{array}
\]
\caption{Results of Algorithm \ref{algoritmo} for the sine transform $F^{(s)}(\omega)$, with $f(x) = \frac{1}{1+e^{\delta x}}$, $ \delta = 1.5$, and $N_1=10$.} \label{table3}
\end{table}

\begin{table} 
\[
\begin{array}{cccccccccccc}
\toprule
\multicolumn{2}{c}{} & \multicolumn{1}{c}{\eta=1e-7} & \multicolumn{1}{c}{\eta=1e-10} & \multicolumn{1}{c}{\eta=1e-13} \\
\midrule
 & {\rm error} & 7.68e-8 & 4.03e-11 & 5.66e-14 \\
\omega = 1 & h & 1.16e-1 & 9.79e-2 & 7.42e-2 \\
& N & 15 & 21 & 31 \\
\midrule
 & {\rm error} & 2.15e-8 & 1.12e-10 & 4.68e-13 \\
\omega = 5 & h & 1.24e-1 & 1.03e-1 & 8.22e-2 \\
& N & 14 & 20 & 28 \\
\midrule
 & {\rm error} & 1.49e-8 & 7.26e-11 & 8.00e-13 \\
\omega = 10 & h & 1.33e-1 & 1.08e-1 & 8.52e-2 \\
& N & 13 & 19 & 27 \\
\bottomrule
\end{array}
\]
\caption{Results of Algorithm \ref{algoritmo} for the sine transform $F^{(s)}(\omega)$, with $f(x) = x^{-\frac{1}{2}}$ and $N_1=10$.} \label{table4}
\end{table}

\begin{figure}
    \begin{center}
    \includegraphics[scale=0.35]{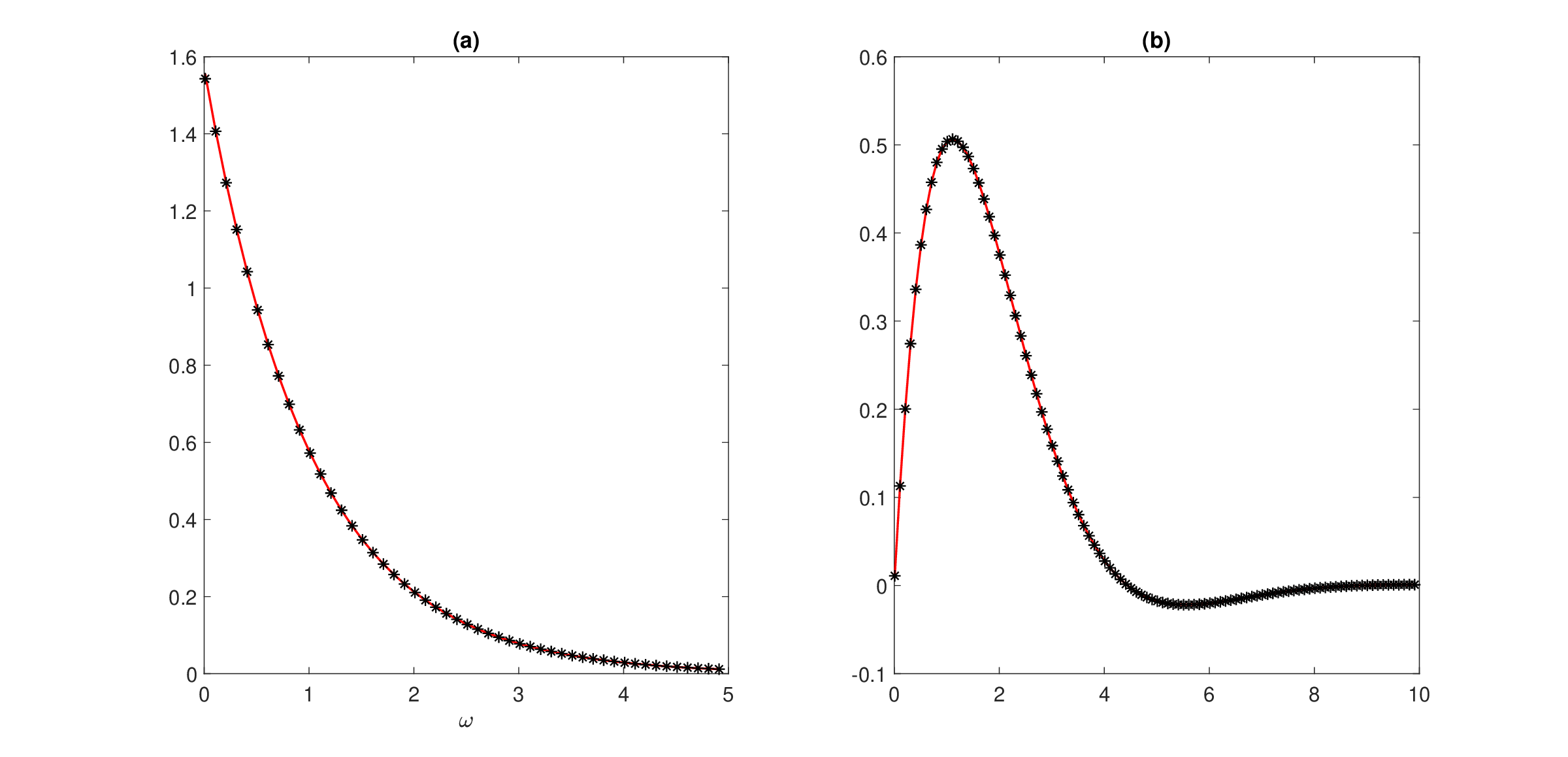}
    \end{center}
    \caption{(a) The cosine transform of the function $f(x)= \frac{1}{1+x^2}$ (red) and (b) the sine transform of the function $f(x)= \frac{x}{1+x^4}$ (red), with the approximations obtained by using the automatic integrator (black). In both cases $\eta = 1e-5$ and $N_1=20$.} \label{figura_trasf1}
\end{figure}

By using Algorithm \ref{algoritmo}, we also reconstruct some known sine and cosine transforms for $\omega \in \left[ \omega_{\rm min}, \omega_{\rm max} \right]$.
Since the error decreases with increasing $\omega$, the idea is to set $N$ and $h$ by applying Algorithm \ref{algoritmo} for $\omega = \omega_{\rm min}$, and to use these fixed values for all the computations.
In this setting, in Figure \ref{figura_trasf1} we plot the cosine transform of the function $f(x) = \frac{1}{1+x^2}$, given by
\begin{equation*}
F^{(c)} \left( \omega \right) = \frac{\pi}{2} e^{-\omega},
\end{equation*}
and the sine transform of $f(x) = \frac{x}{1+x^4}$, that is,
\begin{equation*}
F^{(s)} \left( \omega \right) = \frac{\pi}{2} e^{-\frac{\omega}{\sqrt{2}}} \sin \left( \frac{\omega}{\sqrt{2}} \right),
\end{equation*}
(see \cite[p.408, 3.727, n.1-4]{GR}), together with the approximations computed by employing the automatic integrator.

\section{Conclusion} \label{section 6}

In this work we have considered two double exponential transformations for evaluating Fourier type integrals. In both cases we have derived fairly accurate a priori error estimates, not given in the papers where originally such methods were introduced (\cite{OM,OM1,OM2}).
We have also presented a strategy for automatic integration, together with the corresponding algorithm, which does not require any knowledge of the properties of the function involved in the problem.

\section*{Acknowledgements}

This work was partially supported by GNCS-INdAM and FRA-University of Trieste. The authors are member of the INdAM research group GNCS.


\begin{thebibliography}{99}


\bibitem{AH} A. Asheim and D. Huybrechs, Complex Gaussian quadrature for oscillatory integral transforms, \emph{IMA Journal of Numerical Analysis}, \textbf{33(4) (2013)}, 1322–1341.
 
\bibitem{B} W. Barrett, Convergence of Gaussian quadrature formulae, \emph{Comput. J.} \textbf{3} (1960/1961) 272-277.


\bibitem{DE} J. D. Donaldson and D. Elliott, A unified approach to quadrature rules with asymptotic estimates of their remainders, \emph{SIAM Journal on Numerical Analysis} \textbf{9} (1972) 573-602.

\bibitem{GR} I.S. Gradshteyn and I.M. Ryzhik, \emph{Tables of Integrals, Series, and Products}, 4th ed., Academic Press, New York, 1980.

\bibitem{LB} J. Lund and K. L. Bowers, \emph{Sinc Methods for Quadrature and Differential Equations}, Society for Industrial and Applied Mathematics (SIAM), Philadelphia, 1992.

\bibitem{MS} G., V., Milovanović and M. P. Stanić, \emph{Numerical integration of highly oscillating functions}, Analytic Number Theory, Approximation Theory, and Special Functions: In Honor of Hari M. Srivastava, 2013, 613-649.

\bibitem{NIST} F. Olver, D. Lozier, R. Boisvert and C. Clark, The NIST Handbook of Mathematical Functions, Cambridge University Press, New York, NY (2010).

\bibitem{O} T. Ooura, A Double Exponential Formula for the Fourier Transforms, \emph{Publications of The Research Institute for Mathematical Sciences} \textbf{41} (2005), 971-977.

\bibitem{OM} T. Ooura, M. Mori, The double exponential formula for oscillatory functions over the half infinite interval, \emph{J. Comput. Appl. Math.} \textbf{38} (1991), 353–360.

\bibitem{OM1} T. Ooura, M. Mori, \emph{Double exponential formula for Fourier type integrals with a divergent integrand}, Contributions in Numerical Mathematics, World Scientific Series in Applicable Analysis, World Scientific, Singapore, Vol. 2, 1993, pp. 301–308.

\bibitem{OM2} T. Ooura, M. Mori, A robust double exponential formula for Fourier type integrals, \emph{J. Comput. Appl. Math.} \textbf{112} (1999), 229–241.

\bibitem{Y} P.Yip, \emph{Sine and Cosine Transforms}, in A.D. Poularikas (2nd ed), \emph{The Transforms and Applications Handbook}, Boca Raton, CRC Press LLC, 2000.

\bibitem{Wong} R. S. C. Wong, Quadrature formulas for oscillatory integral transforms, \emph{Numerische Mathematik}, \textbf{39} (1982), 351-360. 


\end{thebibliography}
\end{document}